\documentstyle[amscd,amssymb,verbatim,amsmath,amsthm,amsopn,epsfig,amsfonts,mathrsfs,graphics,pstricks,latexsym,makeidx,12pt]{report}

\theoremstyle{plain}

\theoremstyle{definition}

\numberwithin{equation}{section}
\numberwithin{table}{section}

 \newcommand{\bnum}{\begin{enumerate}}
 \newcommand{\enum}{\end{enumerate}}

\DeclareMathOperator{\Irr}{Irr}

\DeclareMathOperator{\Cl}{C\ell} 

 \DeclareMathOperator{\cd}{cd}

\begin{document}
 \thispagestyle{empty}
 \begin{center}
\bfseries{\LARGE  COMMUTATIVITY DEGREE,}\\     \end{center}
\begin{center}
\bfseries{\LARGE  ITS GENERALIZATIONS,}\\      \end{center}
\begin{center}
\bfseries{\LARGE  AND}\\     \end{center}
\begin{center}\bfseries{\LARGE  CLASSIFICATION OF FINITE   GROUPS} \\   \end{center}
\vspace{.5cm}
\renewcommand{\thefootnote}{}
\footnote{\textit{emails}: rajatkantinath@yahoo.com; akdasnehu@gmail.com}

\begin{center} {\bf\Large (Abstract)} \end{center}
\vspace{.25cm}
\begin{center}
{\large Submitted by}
\end{center}
\vspace{.25cm}
\begin{center}{\bf\Large Rajat Kanti Nath}\\  \end{center}
\vspace{.25cm}
\begin{center}
{\large Under the supervision of} \end{center}
\begin{center} {\bf\Large   Ashish Kumar Das}
\end{center}
\vspace{.25cm}

\begin{center}{\large  In partial fulfilment of the requirement of the degree of}\end{center}

\begin{center} \bfseries {\bf\Large Doctor of Philosophy  in  Mathematics}\\ 
\end{center}

\begin{center}{ \large To}\end{center}

\begin{center}{\bf \large NORTH-EASTERN HILL UNIVERSITY} \end{center}
\begin{center}{\large SHILLONG -- 793022, INDIA}  \end{center}
\begin{center}{\large JULY, 2010 }  \end{center}

 \newpage

Classification of finite groups is a central problem in  theory of groups. Even though finite abelian groups have been completely classified, a lot still remains to be done as far as
non-abelian groups are concerned. People all over the world have used
various types of invariants for classifying finite groups, particularly
the non-abelian ones. The commutativity degree of a finite group is one such invariant, and it seems that many interesting results are possible to obtain with the help of this notion and its generalizations.

In recent years there has been a growing interest in the use of the probabilistic methods in the theory of finite groups. These  methods have proved useful in the solution of several difficult problems on groups. In some cases the probabilistic nature of the problem is apparent from its formulation, but in other cases the connection to probability seems surprising and can not be anticipated by the nature of the problem. 

   The roots of the subject matter of this thesis lie in a series of papers by P. Erd$\rm\ddot o$s and P. Tur$\rm\acute a$n (see \cite{pEpT65, pEpT67, pEpT67b, pEpT68}) published between 1965 and 1968, and also in the Ph. D thesis of K. S. Joseph \cite{ksJ69} submitted in 1969, wherein  some problems of  statistical group theory and commutativity in non-abelian groups have been considered. In 1973, W. H. Gustafson \cite{whG73}  considered the  question -- {\it what is the probability that two group elements commute?} The answer is given by what is known as the commutativity degree of a  group. It may be mentioned here that the question, in some sense, was also considered by Erd$\rm\ddot o$s and  Tur$\rm\acute a$n \cite{pEpT68}. 
   
   Formally, the commutativity degree of a finite group $G$, denoted by $\Pr(G)$,  is defined as the ratio
   
  \[
    \Pr(G) =  \frac{\text{Number of ordered pairs $(x,y) \in G\times G$ such
     that $xy =yx$}}{ \text{Total number of ordered pairs $(x,y) \in G\times
     G$}}.
  \]

\noindent In other words, commutativity degree is a kind of measure for abelianness of a group. Note that  $\Pr(G) > 0$, and that $\Pr(G) = 1$ if and only if $G$ is abelian. Also, given a finite group $G$, we have $\Pr(G) \leq \frac{5}{8}$ with equality if and only if $\frac{G}{Z(G)}$ has order $4$ (see \cite{smBgjS94, whG73}), where $Z(G)$ denotes the center of $G$. This gave rise to the problem of determining the numbers in the interval              $(0, \frac{5}{8}]$ which can be realized as the commutativity degrees of some  finite groups, and also to the problem of classifying  all finite groups with a given commutativity degree.

In 1979, D. J. Rusin \cite{djR79}  computed,  for  a finite group $G$, the values of $\Pr(G)$ when $G' \subseteq Z(G)$, and also when $G' \cap Z(G)$ is trivial,  where $G'$ denotes the commutator subgroup of $G$. He  determined all numbers lying in the interval $(\frac{11}{32}, 1]$ that can be realized as the commutativity degree of some  finite groups, and also classified all finite groups  whose commutativity degrees lie in the interval $(\frac{11}{32}, 1]$.
 
 In 1995, P. Lescot \cite{pL95}   classified, up to isoclinism,  all finite groups whose commutativity degrees are greater than or equal to  $\frac{1}{2}$. It may be mentioned here that the concept of isoclinism between any two groups was introduced by P. Hall \cite{pH40}.  A pair $(\phi,\;\psi )$ is said to be an \textit{isoclinism}  from a  group $G$ to another   group $H$ if the following conditions hold:
\bnum
          \item  $\phi$ is an isomorphism from $G/Z(G)$ to $H/Z(H)$,

          \item  $\psi$ is an isomorphism from $G'$ to $H'$, \; and

          \item the  diagram
\begin{figure}[h]
\setlength{\unitlength}{1mm}
       \begin{picture}(110,42)(-25,15)
   \put(0,46){$\dfrac{G}{Z(G)}  \times \dfrac{G}{Z(G)}$}
   \put(29,47){\vector(1,0){28}}
\put(36,50){$ \phi \times \phi $}
 \put(60,46){$\dfrac{H}{Z(H)}  \times \dfrac{H}{Z(H)}$}
 \put(11,18){$G'$}
   \put(29,19){\vector(1,0){28}}
\put(41,22){$ \psi $}
 \put(71,18){$H'$}
  \put(13,40){\vector(0,-1){16}}
  \put(73,40){\vector(0,-1){16}}
 \put(7,32) {$a_G$}
\put(76,32){$a_H$}
    \end{picture}
     
\end{figure}

\noindent commutes, that is, $ a_H \circ (\phi\times \phi )= \psi\circ
a_G $,   where $a_G$ and  $a_H$  are given respectively by   $a_G  (g_1 Z(G), g_2 Z(G))
= [g_1, g_2]$ for all $g_1, g_2 \in G$ and  $a_H  (h_1 Z(H), h_2 Z(H)) = [h_1, h_2]$ for all $h_1, h_2 \in H$. Here, given $x,y \in G$, $[x,y]$ stands for the commutator $xyx^{-1}y^{-1}$ of $x$ and $y$ in $G$.
\enum

\noindent   In $2001$, Lescot \cite{pL01} has also classified, up to isomorphism, all finite groups whose commutativity degrees  lie   in the interval  $[\frac{1}{2}, 1]$.

 In 2006, F. Barry, D. MacHale and  $\rm\acute A$. N$\rm\acute i$ Sh$\rm\acute e$ \cite{fBdM06}  have shown  that if $G$ is a finite group with    $|G|$   odd and $\Pr(G)>\frac{11}{75}$, then $G$ is 
supersolvable. They also proved that if $\Pr(G) > \frac{1}{3}$, then $G$ is supersolvable. It may be mentioned here that a group $G$ is said to be \textit{supersolvable}  if there is a series of the form
    \[
    \{1\} = A_0\subseteq A_1\subseteq A_2\subseteq\cdots \subseteq A_r = G,
    \]
where $A_i \unlhd G$   and    $A_{i+1}/A_i$ is cyclic for each $i$ with $0 \leq i \leq r-1$.

 In the same year $2006$, R. M. Guralnick and G. R. Robinson  \cite{rmGgrR06} re-established a result of Lescot (see \cite{rmGgrR08}) which says that if $G$ is a finite group with     $\Pr(G)$ greater than $\frac{3}{40}$, then either $G$ is solvable, or  $G \cong A_5 \times B$, where $A_5$ is the alternating group of degree $5$ and $B$ is some abelian group. 

 The classical notion of commutativity degree has been generalized in a number of ways. In 2007,  A. Erfanian, R. Rezaei and P. Lescot \cite{aErRpL07}  studied the probability $\Pr(H, G)$ that an element of a given subgroup $H$ of a finite group $G$ commutes with an element of $G$. Note that $\Pr(G, G)= \Pr(G)$. In 2008, M. R. Pournaki and R. Sobhani \cite{mrPrS08} studied the probability ${\Pr}_g(G)$ that the  commutator of an arbitrarily chosen pair of  elements in a finite group $G$ equals a given group element $g$.  They have also extended some of the results obtained by Rusin. It is easy to see that ${\Pr}_g(G)= \Pr(G)$ if $g=1$, the identity element of $G$.

In Chapter $1$, we briefly recall a few definitions and well-known  results from several relevant topics, which constitute the minimum prerequisites for the subsequent  chapters.  In this chapter, we also fix certain notations.  Given a subgroup $K$ of a group $G$ and an element $x \in G$, we write  $C_K(x)$  and $\Cl_K(x)$ to denote the sets   $\{k \in K : kx =xk \}$ and  $\{k x k^{-1} \in G : k \in K \}$ respectively; noting that, for $K =G$,   these sets coincide respectively with the centralizer  and the conjugacy class of  $x$ in  $G$. Also, given any  two subgroups $H$ and $K$  of a group $G$,  we write $C_K(H) = \{k \in K : hk =kh \text{ for all } h \in H \}$. Note that $C_K(x) = C_K(\langle x \rangle)$, where $\langle x \rangle$ denotes the cyclic subgroup of $G$ generated by $x \in G$.

Further, we write $\Irr (G)$ to denote the set of all irreducible complex characters of $G$, and $\cd(G)$ to denote the set $\{\chi(1) : \chi \in \Irr(G)\}$. If $\chi(1) = |G : Z(G)|^{1/2}$ for some $\chi \in \Irr(G)$, then the group $G$ is said to be of \textit{central type}.

In Chapter $2$, which is based on our papers  \cite{rkNakD12} and \cite{rkNakD10}, we determine, for a finite group $G$, the value of $\Pr(G)$ and the size of $\frac{G}{Z(G)}$ when $|G'| = p^2$ and $|G' \cap Z(G)| = p$, where $p$ is a prime such that $\gcd(p - 1, |G|) = 1$. The main result of Section $2.2$ is given as follows.

\noindent{\bf Theorem 2.2.6.} {\it
Let $G$ be a finite group and $p$ be a prime such that \qquad $\gcd (p-1,|G|) = 1$. If $|G'|= p^2$ and $|G' \cap Z(G)| = p$,  then
\begin{enumerate}
\item $\Pr(G) = 
\begin{cases}
  \frac{2p^2 - 1}{p^4}  &\text{  if $C_G (G')$ is abelian}\\
\frac{1}{p^4} \left( \frac{p-1}{p^{2s-1}} + p^2 +p -1 \right)     &\text{  otherwise,}
\end{cases}$ \smallskip
\item $ |\frac{G}{Z(G)} | =
\begin{cases}
 p^3  &\text{  if $C_G (G')$ is abelian}\\
 p^{2s+2} \text{ or } p^{2s+3}   &\text{  otherwise,}
\end{cases}$\smallskip
\end{enumerate}
where   $p^{2s }=|C_G (G') : Z(C_G (G'))|$.  Moreover, 
\[ \textstyle{
|\frac{G}{G' \cap Z(G)}  : Z(\frac{G}{G' \cap Z(G)})|= |\frac{G}{Z(G)}:Z(\frac{G}{Z(G)})|=p^2.}
\]
}

\noindent This theorem  together with few other supplementary results, enable us to classify all finite groups $G$ of odd order  with $\Pr(G) \geq \frac{11}{75}$. In the process we also   point out a few  small but significant lacunae in the work of Rusin \cite{djR79}. The main result of Section $2.3$ is given as follows.

\noindent{\bf Theorem 2.3.3.} {\it
Let $G$ be a finite group. If $|G|$ is odd and $\Pr(G) \geq \frac{11}{75}$, then the possible values of  $\Pr(G)$ and the corresponding structures of $G'$, $G'\cap Z(G)$ and $G/Z(G)$ are given as  follows:
\begin{table}[h]
\begin{center}
 \begin{tabular}[b]{|c|c|c|c|}
\hline   $\Pr(G)$ & $G'$  & $G'\cap Z(G)$  & $G/Z(G)$ \\ 
\hline $1$ & $\{1 \}$ & $\{1 \}$ & $ \{1 \}$\\  
\hline    $\frac{1}{3}(1+\frac{2}{3^{2s}})$    &  $C_3$ & $ C_3 $ & $ (C_3 \times C_3)^s, \; s \geq 1 $\\   
\hline     $ \frac{1}{5}(1+\frac{4}{5^{2s}}) $    & $C_5$  & $  C_5$ & $ (C_5 \times C_5)^s, \; s \geq 1  $\\    
\hline      $ \frac{5}{21} $	 & $C_7$       & $ \{1 \} $ & $ C_7  \rtimes C_3 $\\    
\hline      $ \frac{55}{343} $	 & $C_7$       & $ C_7 $ & $ C_7 \times C_7$\\   
\hline       $ \frac{17}{81} $	 & $C_9$ or $C_3 \times C_3$     & $ C_3  $ &    $(C_3 \times C_3)  \rtimes C_3 $\\   
&$C_3 \times C_3 $ & $ C_3 \times C_3$ & $C_3 ^3$\\   
\hline     $ \frac{121}{729} $	 & $C_3 \times C_3$       & $ C_3 \times C_3 $ &    $ C_3 ^4 $\\   
\hline      $ \frac{7}{39} $ 	 & $C_{13}$       &$ \{1 \} $  &   $ C_{13} \rtimes C_3 $ \\   
\hline      $ \frac{3}{19} $	 &  $C_{19}$      & $ \{1 \} $ &   $C_{19} \rtimes C_3$ \\   
\hline      $ \frac{29}{189} $ 	 & $C_{21}$       & $ C_3 $ &   $C_3\times ( C_7 \rtimes C_3)  $ \\   
\hline      $ \frac{11}{75} $	 & $C_5 \times C_5$       & $ \{1 \} $ &    $ (C_5 \times C_5)  \rtimes C_3$ \\   
\hline
\end{tabular}
\end{center} 
 
\end{table}}
\newpage
\noindent In the above table  $C_n$ denotes the cyclic group of order $n$ and $\rtimes$ stands for semidirect product.

In \cite[Corollary 2.3]{mrPrS08},  M. R. Pournaki and R. Sobhani have proved that, for a finite group $G$ satisfying $|\cd(G)|=2$, one has
\[
 \Pr(G) \geq \dfrac{1}{|G'|} \left( 1 + \dfrac{|G'| - 1}{|G : Z(G)|}\right) 
\]
with equality   if and only if $G$ is of central type. In Section $2.4$, we have improved this result as follows.
  
\noindent{\bf Theorem 2.4.1.} { \it
If $G$ is a finite group, then
\[
 \Pr(G) \geq \dfrac{1}{|G'|} \left( 1 + \dfrac{|G'| - 1}{|G : Z(G)|}\right).
\]
In particular,  $\Pr(G) > \frac{1}{|G'|}$ if $G$ is non-abelian.}

There are  several equivalent conditions that are necessary as well as sufficient for the attainment of the above lower bound for $\Pr(G)$.
 
\noindent{\bf Theorem 2.4.3.} {\it
For a finite non-abelian group $G$, the statements given below are equivalent.
\begin{enumerate}
   \item $\Pr(G) = \dfrac{1}{|G'|} \left( 1 + \dfrac{|G'| - 1}{|G : Z(G)|}\right)$.     
\item  $\cd(G) = \{1, |G : Z(G)|^{1/2} \}$, which means that  $G$ is of central type with $|\cd(G)| = 2$.      
\item $|\Cl_G(x)|  = |G'|$ for all $x \in G - Z(G)$.    
\item $\Cl_G(x)  = G'x$ for all $x \in G - Z(G)$; in particular, $G$ is a nilpotent group of class $2$.    
\item $C_G(x) \unlhd G$ and $G' \cong \frac{G}{C_G(x)}$ for all $x \in G - Z(G)$; in particular, $G$ is a $CN$-group, that is, the centralizer of every element is normal.      
\item $G' = \{[y,x]  : y \in G\}$ for all $x \in G - Z(G)$; in particular, every element of $G'$ is a commutator.  
\end{enumerate}
}
\noindent  Theorem $2.4.1$ and Theorem $2.4.3$ not only allow us to obtain some characterizations for finite nilpotent groups of class $2$ whose commutator subgroups have prime order, but also enable us to re-establish certain facts (essentially due to K. S. Joseph \cite{ksJ69})  concerning    the smallest prime divisors of the orders of finite groups.

\noindent{\bf Proposition 2.4.4.} {\it
Let $G$ be a finite group and $p$ be the smallest prime divisor of $|G|$.
\begin{enumerate}
\item If $p \ne 2$, then $\Pr(G) \neq \frac{1}{p}$.    
\item When $G$ is non-abelian, $\Pr(G) > \frac{1}{p}$ \, if and only if \, $|G'| = p$ and $G' \subseteq Z(G)$.   
\end{enumerate}
}

\noindent{\bf Corrolary 2.4.5.} {\it
If $G$ is a finite group with $\Pr(G) = \frac{1}{3}$, then $|G|$ is even.
}

\noindent{\bf Proposition 2.4.7.} {\it
Let $G$ be a finite group and $p$ be a prime. Then the following statements are equivalent.
\begin{enumerate}
\item $|G'| = p$ and $G' \subseteq Z(G)$.    
\item $G$ is of central type with $|\cd(G)| = 2$ and $|G'| = p$.  
\item G is a direct  product of a $p$-group $P$ and an abelian group $A$ such that $|P'| = p$ and $\gcd (p, |A|) =1$.   
\item $G$ is isoclinic to an extra-special $p$-group; consequently, $|G:Z(G)| = p^{2k}$ for some positive integer  $k$.  
\end{enumerate}
In particular, if $G$ is non-abelian and $p$ is the smallest prime divisor of $|G|$, then  the above statements are also equivalent to the condition $\Pr(G) > \frac{1}{p}$.
}

In \cite{pL01}, Lescot deduced that  $\Pr(D_{2n}) \rightarrow \frac{1}{4}$ and $\Pr(Q_{2^{n +1}}) \rightarrow \frac{1}{4}$ as $n \rightarrow \infty$, where $D_{2n}$ and $Q_{2^{n +1}}$  denote the dihedral group of  order $2n$, $n \geq 1$, and the quaternion group of order $2^{n+1}$, $n \geq 2$, respectively.  He also enquired whether there are other natural families of finite groups with the same property. In 2007,  I. V. Erovenko and B. Sury \cite{ivEbS08}  have shown, in particular, that for every integer $k > 1$ there exists a family $\{G_n\}$ of finite groups such   that  $\Pr(G_n) \rightarrow \frac{1}{k^2}$ as $n  \rightarrow \infty$. In the last section of Chapter $2$ we have considered the question posed by Lescot mentioned above. Moreover, we make the following observation.

\noindent{\bf Proposition 2.5.1.} {\it
For every integer $k > 1$ there exists a family $\{G_n\}$ of finite groups  such that $\Pr(G_n) \rightarrow \frac{1}{k}$ as $n  \rightarrow \infty$.}

\noindent In the same line, we also have the following result.
  
\noindent{\bf Proposition 2.5.2.} {\it  
For every positive integer $n$ there exists a finite group $G$ such that $\Pr(G) = \frac{1}{n}$.
}
  
In Chapter $3$, which is based on our papers \cite{akDrkN09} and  \cite{akDrkN10},  we generalize the following result of F. G. Frobenius \cite{fgF68}:

   \textit{If $G$ is a finite group and  $g \in G$, then the number of solutions of the commutator equation $xyx^{-1}y^{-1} = g$ in $G$ defines a character on $G$, and is given by
\[
    \zeta (g) = \underset{\chi\in\Irr(G)}{\sum}\frac{|G|}{\chi(1)}\chi(g).
\]}

  We write \, $F(x_1, x_2, \dots, x_n)$ \, to denote the free group of words  on $n$ \quad generators  $x_1, x_2, \dots, x_n$.  For $1\leq i\leq n$, we  write  `$x_i\in \omega(x_1, x_2, \dots, x_n)$'  to mean that $x_i$ has a non-zero index (that is, $x_i^k$ forms a syllable, with $0 \ne k  \in \mathbb{Z}$) in the word $\omega(x_1, x_2, \dots, x_n) \in F(x_1, x_2, \dots, x_n)$. We call a word $\omega(x_1, x_2, \dots, x_n)$     \textit{admissible} if each $x_i\in\omega(x_1, x_2, \dots, x_n)$
has precisely two non-zero indices, namely, $+1$ and $-1$.  We write ${\mathscr{A}}(x_1, x_2, \dots, x_n)$  to denote the set of all admissible words in $F(x_1, x_2, \dots, x_n)$. 

Given a finite group $G$ and an element $g \in G$, let
${\zeta}^{\omega}_n(g)$ denote the number of solutions
$(g_1, g_2, \dots, g_n)\in G^n$ \, of the word equation  $\omega(x_1, x_2, \dots, x_n)=g$, where  $G^n = G \times G \times \dots \times G$ ($n$ times).  Thus,
\[
   {\zeta}^{\omega}_n(g) = |\{(g_1, g_2, \dots, g_n)\in G^n :
    \omega(g_1, g_2, \dots, g_n) = g\}|.
\]
The main result of Section $3.1$ is given as follows.

\noindent{\bf Theorem 3.1.4.} {\it
Let  ${\omega}(x_1, x_2, \dots, x_n)\in {\mathscr{A}}(x_1, x_2, \dots,x_n)$, $n \geq 1$. If $G$ is a finite group, then   the map ${\zeta}^{\omega}_n : G \longrightarrow
{\mathbb{C}}$  defined by
\[
    {\zeta}^{\omega}_n(g) = |\{(g_1, g_2, \dots, g_n)\in G^n:
      {\omega}(g_1, g_2, \dots, g_n)=g\}|,\;\; g\in G,
\]
is a character of $G$.
}

  Given any two finite sets $X$ and $Y$, a function $f : X \longrightarrow Y$ is said to be \textit{almost measure preserving} if there exists a sufficiently small positive
real number $\epsilon$ such that
\[
\left| \dfrac{|f^{-1}(Y_0)|}{|X|}  - \dfrac{|Y_0|}{|Y|}\right| < \epsilon   \quad \text{for all} \; Y_0 \subseteq Y.
\] 

In Section $3.2$, we consider the following question posed by Aner Shalev \cite[Problem 2.10]{aS07}:

\noindent \textit{Which words induce almost measure preserving maps on finite simple groups?} 

\noindent More precisely, given an admissible word $\omega(x_1, x_2, \dots, x_n)$ and the induced word map $\alpha_{\omega} : G^n \to G$ defined by 
$\alpha_{\omega} (g_1, g_2, \dots, g_n) = \omega (g_1, g_2, \dots, g_n)$, we proved that

\noindent{\bf Corollary 3.2.7.} {\it
Let $G$ be a finite simple group, and $o(1)$ be a real number depending on $G$ which tends to zero as $|G| \rightarrow \infty$.  
\bnum
   \item If \, $Y \subseteq G$, \, then \, $\dfrac{|(\alpha_{\omega})^{- 1}(Y)|}{|G|^n} = \dfrac{|Y|}{|G|}  + o(1)$.   This means that
   the  map $\alpha_{\omega}$ is almost measure preserving. 
   \item If $X \subseteq G^n$, then
   $\dfrac{|\alpha_{\omega}(X)|}{|G|} \geq \dfrac{|X|}{|G|^n}  - o(1)$; in particular, 
 if $X$ is such that $|X| = (1 - o(1))|G|^n$,  then  $ |\alpha_{\omega}(X)| =  (1 - o(1))|G|$.  This means that almost all the elements of $G$ can be expressed as $\omega(g_1, g_2, \dots, g_n)$ for some $g_1, g_2, \dots, g_n \in G$.
\enum
}

In the last section of Chapter $3$, we obtain yet another generalization of Frobenius' result mentioned above. The main results of this section are given as follows.

\noindent{\bf Theorem 3.3.1.} {\it
Let $G$ be a finite group, $H \unlhd G$ and $g \in G$. If  $\tilde{\zeta}(g)$ denotes the number of elements
$(h_1, g_2) \in H \times G$ satisfying  $[h_1, g_2] = g$, then 
$\tilde{\zeta}$ is a class function of $G$ and 
\[
 \tilde{\zeta}(g) = \underset{\chi\in \Irr(G)}{\sum} \dfrac{|H| [{\chi}_{_H},{\chi}_{_H}]}{\chi(1)}\chi(g) = \underset{\chi\in \Irr(G)}{\sum} \dfrac{|H| [{\chi}_{_H}^G,{\chi}]}{\chi(1)}\chi(g).
\]
}

\noindent{\bf Corollary 3.3.2.} {\it
Let $G$ be a finite group. Then, with notations as above,  $\tilde{\zeta}$ is a character of $G$.
}

\noindent{\bf Proposition 3.3.3.} {\it
Let $G$ be a finite group, $H \unlhd G$ and $g \in G$. If \;  $\tilde{\zeta}_{2n}(g)$, $n \geq 1$, denotes the number of elements $((h_1,g_1),  \dots, (h_n,g_n)) \in (H  \times G)^n$ satisfying  $[h_1, g_1] \dots [h_n, g_n] = g$,  then  $\tilde{\zeta}_{2n}$ is a character of $G$ and 
\[
\tilde{\zeta}_{2n}(g) = \underset{\chi\in \Irr(G)}{\sum} \dfrac{|G|^{n - 1}|H|^n[{\chi}_{_H},{\chi}_{_H}]^n}{\chi(1)^{2n - 1}} \chi(g).
\]
}

\noindent\textbf{Proposition 3.3.4.} {\it
Let $H$ be a  subgroup of a finite group $G$ and $g \in G$. \!Then the number of elements $(g_1, h_2, g_3) \in G \times H \times G$ satisfying
$g_1h_2g_1^{-1}g_3h_2^{-1}g_3^{-1} = g$  defines a character of $G$ and is given by
\[
\tilde{\zeta}_3(g) = \underset{\chi\in \Irr(G)}{\sum} \dfrac{|G||H|[{\chi}_{_H},{\chi}_{_H}]}{\chi(1)} \chi(g).
\]
}

  In Chapter $4$, which is based on our paper \cite{rkNakD11},  we study the probability ${\Pr}_g^{\omega}(G)$  that an arbitrarily chosen $n$-tuple of elements of a given finite group $G$ is mapped to a given group element $g$ under the word map induced by a non-trivial admissible word $\omega (x_1, x_2, \dots, x_n)$. Formally, we write
\begin{equation*} 
{\Pr}_g^{\omega}(G) = \dfrac{{\zeta}^{\omega}_n (g)}{|G^n|},
\end{equation*}
where ${\zeta}^{\omega}_n (g) =  |\{(g_1, g_2, \dots, g_n) \in G^n : \omega (g_1, g_2, \dots, g_n) = g \}|$.   The main results of Section $4.1$ are as follows.
  
\noindent{\bf Proposition 4.1.1.} {\it   
Let $G$ be a finite group and $\omega(x_1, x_2, \dots, x_n)$  be a non-trivial admissible word. Then
\bnum
\item  \  ${\Pr}^{\omega}_1(G)\;  \geq \;  \dfrac{n|G:Z(G)|- n +1}{|G:Z(G)|^n} \; \geq \; \dfrac{1}{|G:Z(G)|^n} \; \gneq \;  0$, \label{pr1nz} \vspace{0.25cm}
\item  \  ${\Pr}^{\omega}_1(G) = 1$ \text{  if and only if   } $G$ \text{  is abelian}. \label{gabel}
\enum
}

\noindent{\bf Proposition 4.1.3.} {\it
Let $G$ and $H$ be two finite groups and    $(\phi, \psi)$ be an isoclinism from $G$ to $H$. If $g \in G'$ and $\omega(x_1, x_2, \dots, x_n)$  is a non-trivial admissible word,    then  
 \[
 {\Pr}_g^{\omega}(G) = {\Pr}_{\psi(g)}^{\omega}(H). 
 \]
}  

\noindent{\bf Proposition 4.1.4.} {\it
Let $G$ be a finite group and $\omega(x_1, x_2, \dots, x_n)$  be a non-trivial admissible word.  If $g, h \in G'$   generate the same cyclic subgroup of $G$, then  
${\Pr}^{\omega}_g(G) = {\Pr}_h^{\omega}(G)$.
}  

\noindent{\bf Proposition 4.1.6.} {\it
Let $G$ be a finite group, $g \in G'$ and $\omega(x_1, x_2, \dots, x_n)$  be a non-trivial admissible word.  Then
\bnum
\item    ${\Pr}^{\omega}_g(G) \leq {\Pr}^{\omega}_1(G) \leq {\Pr}(G)$, \label{prwbd1}
\item   ${\Pr}^{\omega}_g(G) = {\Pr}^{\omega}_1(G) $ if and only if $g = 1$.  \label{prwbd2}
\enum
}
  
 Let   $m_G = \min \{ \chi (1) \, : \, \chi \in \Irr (G), \chi (1) \ne 1 \}$.  Considering, in particular,    the word $x_1 x_2\dots x_n x_1^{-1} x_2^{-1}\dots x_n^{-1}$, $n \geq 2$, and writing ${\Pr}^n_g(G)$ in place of ${\Pr}^{\omega}_g(G)$,  we obtain the following results in the 
sections $4.2$,  $4.3$ and    
 $4.4$.

\noindent{\bf Proposition 4.2.2.} {\it
Let $G$ be a finite non-abelian group, $g \in G'$ and  $d$  be an integer such that $2 \leq d \leq m_G$.  Then  
\bnum
\item  \  $\left|{\Pr}^{n}_g(G) - \dfrac{1}{|G'|}\right|    \leq \dfrac{1}{d^{n-2}} \left( \Pr(G) - \dfrac{1}{|G'|}\right)$. \quad In other words, \\

  \  $\dfrac{1}{d^{n-2}}\left(
-\Pr(G) + \dfrac{d^{n-2} + 1}{|G'|}\right) \leq  {\Pr}^{n}_g(G) \leq \dfrac{1}{d^{n-2}} \left( \Pr(G) + \dfrac{d^{n-2} - 1}{|G'|}\right)$. \\
 
\item   \  $\left|{\Pr}^{n}_g(G) - \dfrac{1}{|G'|}\right|  \leq  \dfrac{1}{d^n} \left( 1 - \dfrac{1}{|G'|}\right)$.  \quad In other words, \\

  \  $\dfrac{1}{d^n}\left( -1 + \dfrac{d^n + 1}{|G'|}\right) \leq  {\Pr}^{n}_g(G)  \leq  \dfrac{1}{d^n}\left(1 + \dfrac{d^n -1}{|G'|}\right) $. \\
  
  \  In particular, ${\Pr}^{n}_g(G)  \leq  \dfrac{2^n +1}{2^{n+1}}$.
\enum
}  

\noindent{\bf Proposition 4.2.3.} {\it
If $G$ is a finite non-abelian simple group and $g \in G'$, then 
\[
\left|{\Pr}^{n}_g(G) - \dfrac{1}{|G|}\right|    \leq \dfrac{1}{3^{n-2}} \left( \dfrac{1}{12} - \dfrac{1}{|G|}\right).
\]
In other words,
\[
\dfrac{1}{3^{n-2}} \left( \dfrac{-1}{12} + \dfrac{3^{n-2} +1}{|G|}\right) \leq  {\Pr}^{n}_g(G)  \leq \dfrac{1}{3^{n-2}} \left( \dfrac{1}{12} + \dfrac{3^{n-2} -1}{|G|}\right).
\]
In particular,
\[
{\Pr}^{n}_g(G) \leq \dfrac{3^{n-2} + 4}{3^{n-1}\times 20}.
\]
}

\noindent{\bf Corollary 4.3.2.} {\it
Let $G$ be a finite non-abelian  group with $|\cd(G)| =2$. Then  every element of $G'$ is a generalized commutator of length $n$ for all $n \geq 2$; in particular, every element of $G'$ is a  commutator.
}

\noindent{\bf Proposition 4.3.3.} {\it
Let $G$ be a finite non-abelian group, $g \in G'$ and  $d$  be an integer such that $2 \leq d \leq m_G$.   Then
\bnum
\item \  $ {\Pr}^{n}_g(G)  \; = \; \dfrac{1}{d^{n-2}} \left( \Pr(G) + \dfrac{d^{n-2} - 1}{|G'|}\right)$ \quad   if and only if  
 
 \  $g=1$ and $\cd(G) = \{1, d\}$.  
 
\item  \  ${\Pr}^{n}_g(G) = \dfrac{1}{d^{n-2}}\left(
-\Pr(G) + \dfrac{d^{n-2} + 1}{|G'|}\right)$ \quad   if and only if  
 
 \  $g \ne 1$, $\cd(G)=\{1, d\}$   and  $|G'|=2$. 
\enum

}

\noindent{\bf Proposition 4.3.4.} {\it
Let $G$ be a finite non-abelian group, $g \in G'$ and  $d$  be an integer such that $2 \leq d \leq m_G$.   Then
\bnum
\item \  $ {\Pr}^{n}_g(G)  \; = \; \dfrac{1}{d^n}\left(1 + \dfrac{d^n -1}{|G'|}\right)$ \quad   if and only if  
 
 \  $g=1$ and $\cd(G) = \{1, d\}$.   
 
\item  \ ${\Pr}_g^n(G) \; = \; \dfrac{1}{d^n}\left(
-1 + \dfrac{d^n + 1}{|G'|}\right)$  \quad   if and only if  
 
 \  $g \ne 1$, $\cd(G)=\{1, d\}$    and   $|G'|=2$.  
\enum
}

\noindent{\bf Proposition 4.3.7.} {\it
Let $G$ be a finite non-abelian group, $g \in G'$ and  $p$ be the smallest prime divisor of $|G|$. Then
\[   
{\Pr}^{n}_g(G) = \dfrac{p^n + p - 1}{p^{n+1}} 
\]
if and only if \, $g=1$,  and   $G$ is isoclinic to 
\[
\langle x, y  \; : \; x^{p^2} = 1 =y^p, \; y^{-1}xy = x^{p+1}\rangle.
\]
In particular, putting $p=2$, ${\Pr}^{n}_g(G) = \dfrac{2^n + 1}{2^{n+1}}$ \; if and only if \, $g=1$,  and    $G$ is isoclinic to $D_8$, the dihedral group, and hence,  to  $Q_8$, the group  of quaternions.
}

\noindent{\bf Proposition 4.4.1.} {\it
 Let $G$ be a finite non-abelian  group with $|\cd(G)| = 2$ and $g \in G'$.  Then
 \begin{align*}
 {\Pr}^{n}_1(G) &\ge \dfrac{1}{|G'|}\left(1 + \dfrac{|G'| -1}{|G:Z(G)|^{n/2}}\right) \quad  and \\
 {\Pr}^{n}_g(G) &\le \dfrac{1}{|G'|}\left(1 - \dfrac{1}{|G:Z(G)|^{n/2}}\right)\, \quad if \; g\ne 1.
 \end{align*}
 Moreover, in each case, the equality holds if and only if $G$ is of central type.
} 
 
\noindent{\bf Corollary 4.4.2.} {\it
Let $G$ be a finite non-abelian group  and $g \in G'$. If $G$ is of central type with $|\cd(G)|=2$, then
\begin{align*}
 {\Pr}^{n}_1(G) &\le \dfrac{1}{|G'|}\left(1 + \dfrac{|G'| -1}{2^n}\right) \quad  and \\
 {\Pr}^{n}_g(G) &\ge \dfrac{1}{|G'|}\left(1 - \dfrac{1}{2^n}\right)\, \quad if \; g\ne 1.
 \end{align*}
} 
 
\noindent{\bf Proposition 4.4.3.} {\it
Let $G$ be a finite non-abelian group  and $g \in G'$.  If $G' \subseteq Z(G)$ and $|G'| = p$, where  $p$ is a prime, then
\[
        {\Pr}^{n}_g(G)=
        \begin{cases}
             
             \frac{1}{p}\left ( 1+\frac{p-1}{p^{nk}}\right )
                        &\;\text{if \, $g=1$}\\
             \frac{1}{p}\left ( 1-\frac{1}{p^{nk}}\right )
                        &\; \text{if \, $g\ne 1$,} 
        \end{cases}
      \]
where $k= \dfrac{1}{2}\log_p |G:Z(G)|$. 
}

\noindent{\bf Proposition 4.4.5.} {\it
   Let $p$ be  a prime. Let $r$ and $s$ be two positive integers such that $s\mid (p - 1)$,\, and \; $r^j\equiv 1 \pmod p$ \, if and only if \, $s\mid j$.  If \;  $G = \langle a,b : a^p = b^s = 1, bab^{-1} = a^r\rangle $ \; and \; $g\in G'$, then  
\[
{\Pr}^n_g(G)=
        \begin{cases}
             \dfrac{s^n + p - 1}{ps^n}
                    &\; \text{ if   $g=1$}
                    \vspace{.25cm} \\

             \dfrac{s^n - 1}{ps^n}
                    &\; \text{ if   $g\ne 1$}.             
        \end{cases}
\]
}  

\noindent{\bf Proposition 4.4.6.} {\it
 Let $G$ be a finite non-abelian group  and $g \in G'$.  If $G'\cap Z(G)=\{1\}$ and $|G'| = p$, where  $p$ is a prime, then  
\bnum
\item   \quad    $G$ is isoclinic to the group \, $\langle a,b\; :\; a^p = b^s = 1, \, bab^{-1} = a^r\rangle$, \, where  

\quad     $s\mid (p - 1)$,  and $r^j\equiv 1 \pmod p$ if and only if $s\mid j$,   \vspace{.25cm}

\item    \quad   ${\Pr}^n_g(G)=
        \begin{cases}
             \dfrac{s^n + p - 1}{ps^n}
                    &\; \text{ if   $g=1$} \vspace{.25cm} \\

             \dfrac{s^n - 1}{ps^n}
                    &\; \text{ if   $g\ne 1$}.             
        \end{cases}$  
\enum
}  

\noindent{\bf Proposition 4.4.7.} {\it
Let $G$ be a finite non-abelian group and $g \in G'$. If $g \ne 1$, then ${\Pr}_g^n(G) <\frac{1}{p}$, where $p$ is the smallest prime divisor of $|G|$. In particular, we have ${\Pr}_g^n(G) <\frac{1}{2}$.
}  

\noindent{\bf Proposition 4.4.8.} {\it
 For each $\varepsilon > 0$ and for each prime $p$, there exists a finite group $G$ such that  
 \[
 \left|{\Pr}_g^n(G)- \dfrac{1}{p}\right| < \varepsilon
 \]
 for all $g \in G'$.
}

Let $G$ be a finite group and $g \in G'$. Let $H$ and   $K$ be two subgroups of $G$.
 In Chapter $5$, which is based on our paper \cite{akDrkN10},  we study the probability ${\Pr}_g(H, K)$ that  the commutator of a randomly chosen pair of elements (one from $H$ and the other from $K$) equals $g$.  In other words, we study the ratio
\[
{\Pr}_g(H, K) = \frac{|\{(x, y)\in H \times K :  [x,y]  = g\}|}{|H||K|},
 \]  
and further extend some of the results obtained in \cite{aErRpL07} and \cite{mrPrS08}.  Without any loss, we may assume that $G$ is non-abelian. The main results of the  sections $5.1$ and $5.2$  are as follows.
 
\noindent{\bf Proposition 5.1.1.} {\it
Let $G$ be a finite group and $g \in G'$. If $H$ and $K$ are two subgroups of $G$, then
${\Pr}_g(H,K)= {\Pr}_{g^{-1}}(K, H)$. However, if $g^2 =1$, or if $g \in H \cup K$ (for example, when $H$ or $K$ is normal in $G$), we have ${\Pr}_g(H,K) = {\Pr}_g(K, H) = {\Pr}_{g^{-1}}(H, K)$. 
}

\noindent{\bf Theorem 5.1.3.} {\it
Let $G$ be a finite group and $g \in G'$. If $H$ and $K$ are two subgroups of $G$, then
\[
{\Pr}_g (H,K) = \dfrac{1}{|H||K|} \underset{g^{-1}x \in \Cl_K (x)}{\underset{x\in H}{\sum}} |C_K(x)| = \dfrac{1}{|H|} \underset{g^{-1}x \in \Cl_K (x)}{\underset{x\in H}{\sum}} \dfrac{1}{|\Cl_K (x)|}, 
\]
where $C_K(x) = \{y \in K : xy =yx \}$ and $\Cl_K (x) = \{ yxy^{-1} : y \in K \}$, the $K$-conjugacy class of $x$. 
}

\noindent This theorem plays a key role in the study of ${\Pr}_g(H,K)$.
As an immediate consequence, we have the following generalization of the well-known formula  ${\Pr}(G) = \frac{k(G)}{|G|}$. 

\noindent{\bf Corollary 5.1.4.} {\it
Let $G$ be a finite group and  $H$, $K$ be two subgroups of $G$. If $H \unlhd G$, then
\[
{\Pr}(H,K) = \dfrac{k_K (H)}{|H|},
\]
where $k_K(H)$ is the number of $K$-conjugacy classes  that constitute $H$.  
}

\noindent{\bf Proposition 5.1.5.} {\it
If $H$ is an abelian normal subgroup of a finite group $G$ with a complement $K$ in $G$ and $g \in G'$, then
\[
{\Pr}_g(H, G) = {\Pr}_g(H, K).
\]
}

\noindent{\bf Corollary 5.1.6.} {\it
Let $G$ be a finite group and $g \in G'$. If $H \unlhd G$ with $C_G (x) = H$ for all $x \in H - \{1\}$, then 
\[
{\Pr}_g (H,G)= {\Pr}_g (H,K), 
\] 
where $K$ is a complement of $H$ in $G$. In particular,
\[
{\Pr}  (H,G)= \dfrac{1}{|H|} + \dfrac{|H|-1}{|G|}.
\] 
}

\noindent{\bf Proposition 5.2.1.} {\it
Let $G$ be a finite group and $g \in G'$.  Let $H$ and $K$ be any two subgroups of $G$.  If $g \neq 1$, then
\begin{enumerate}
\item ${\Pr}_g(H,K) \neq 0  \Longrightarrow {\Pr}_g(H,K) \geq \dfrac{|C_H(K)||C_K(H)|}{|H||K|}$,   
\item ${\Pr}_g(H,G) \neq 0  \Longrightarrow {\Pr}_g(H,G) \geq \dfrac{2|H \cap Z(G)||C_G(H)|}{|H||G|}$,   
\item ${\Pr}_g(G) \neq 0  \Longrightarrow {\Pr}_g(G) \geq \dfrac{3}{|G: Z(G)|^2}$. 
\end{enumerate}
}

\noindent{\bf Proposition 5.2.2.} {\it
Let $G$ be a finite group and $g \in G'$. If $H$ and $K$ are any two subgroups of $G$, then
\[
{\Pr}_g(H,K) \leq {\Pr}(H,K)
\]
with equality  if and only if $g = 1$.
}

\noindent{\bf Proposition 5.2.3.} {\it
Let $G$ be a finite group and $g \in G'$, $g \neq 1$. Let  $H$ and $K$ be any two subgroups of $G$.  If $p$ is the smallest prime divisor of $|G|$, then
\[
{\Pr}_g(H,K) \leq \dfrac{|H| - |C_H(K)|}{p|H|} < \dfrac{1}{p}.
\]
}

\noindent{\bf Proposition 5.2.4.} {\it
Let $H$, $K_1$ and $K_2$  be subgroups of a finite group $G$ with $K_1 \subseteq K_2$. Then 
\[
{\Pr}(H, K_1) \geq {\Pr}(H, K_2) 
\]
with equality if and only if \; 
$\Cl_{K_1} (x) = \Cl_{K_2} (x)$ \, for all $x \in H$. 
}

\noindent{\bf Proposition 5.2.5.} {\it
Let $H$, $K_1$ and $K_2$  be subgroups of a finite group $G$ with $K_1 \subseteq K_2$.  Then
\[
\Pr(H, K_2) \geq \dfrac{1}{|K_2 : K_1|} \left(\Pr(H, K_1) + \dfrac{|K_2| - |K_1|}{|H||K_1|}\right)
\]
with equality if and only if \;  $C_H(x) = \{1\}$ \,  for all $x \in K_2 - K_1$.
}

\noindent{\bf Proposition 5.2.6.} {\it
Let $H_1 \subseteq H_2$ and $K_1 \subseteq K_2$ be subgroups of a finite group $G$ and $g \in G'$.  Then
\[
{\Pr}_g (H_1,K_1)\leq |H_2:H_1||K_2:K_1|{\Pr}_g (H_2, K_2)
\] 
with equality if and only if  
\begin{align*}
&g^{-1}x \notin \Cl_{K_2} (x) \; \text{ for all } \; x \in H_2 - H_1,\\ &g^{-1}x \notin \Cl_{K_2} (x)-\Cl_{K_1} (x)\; \text{ for all } \; x \in  H_1,\\
 \text{and   } \quad &C_{K_1}(x)=C_{K_2}(x) \; \text{ for all } \; x \in  H_1 \; \text{ with } \; g^{-1}x \in \Cl_{K_1} (x). 
\end{align*}
In particular, for $g=1$, the condition for equality reduces to
 $H_1 = H_2$, and $K_1 = K_2$.
}

\noindent{\bf Corollary 5.2.7.} {\it 
Let $G$ be a finite group, $H$ be a subgroup of $G$ and  $g \in G'$. Then
 \[
 {\Pr}_g (H,G)\leq |G : H| \Pr(G)
\]
with equality if and only if $g=1$ and $H=G$.
}

\noindent{\bf Theorem 5.2.8.} {\it
Let $G$ be a finite group and $p$ be the smallest prime dividing $|G|$. If $H$ and $K$ are any two subgroups of $G$, then
\begin{align*}
{\Pr}(H, K) &\geq \dfrac{|C_H(K)|}{|H|} + \dfrac{p(|H| - |X_H| - |C_H(K)|)+ |X_H|}{|H||K|}    \\
\text{and  } \; {\Pr}(H, K) &\leq \dfrac{(p-1)|C_H(K)| + |H|}{p|H|} - \dfrac{|X_H|(|K| - p)}{p|H||K|}, 
\end{align*}
where $X_H = \{ x \in H : C_K(x) = 1\}$. Moreover, in each of these bounds, $H$ and $K$ can be interchanged.
}

\noindent{\bf Corollary 5.2.9.} {\it
Let $G$ be a finite group and $p$ be the smallest prime dividing $|G|$. If $H$ and $K$ are   two subgroups of $G$ such that $[H, K] \neq \{ 1 \}$, then 
\[
{\Pr}(H, K) \leq \dfrac{2p - 1}{p^2}.
\]
In particular,  ${\Pr}(H, K) \leq \frac{3}{4}$.
}

\noindent{\bf Proposition 5.2.10.} {\it
Let $G$ be a finite group and   $H$, $K$ be any two subgroups of $G$. If ${\Pr}(H, K) = \frac{2p-1}{p^2}$ for some prime $p$, then $p$ divides $|G|$. If $p$ happens to be the smallest prime divisor of $|G|$, then  
 \[
 \dfrac{H}{C_H(K)} \cong C_p \cong \dfrac{K}{C_K(H)}, \, \text{ and hence, } \, H \neq K. 
 \]
In particular,   $\frac{H}{C_H(K)} \cong C_2 \cong \frac{K}{C_K(H)}$ if ${\Pr}(H, K) = \frac{3}{4}$.
}

In the last section of chapter $5$, with $H$ normal in $G$, we also develop and study a character theoretic formula for ${\Pr}_g(H, G)$ given by 
\[
{\Pr}_g (H,G)  = \dfrac{1}{|G|} \underset{\chi\in \Irr(G)}{\sum}
          \dfrac{[ {\chi}_{_H},{\chi}_{_H} ]}{\chi(1)}\chi(g).
\]

\noindent{\bf Proposition 5.3.1.} {\it
Let $G$ be a finite group. If $H$ is a normal subgroup of $G$ and $g \in G'$,   then
\[\left|{\Pr}_g(H, G) - \dfrac{1}{|G'|}\right| \leq |G:H| \left( \Pr(G) - \dfrac{1}{|G'|}\right).
\]
}

As an application, we obtain yet another condition under which every element of   $G'$   is a commutator.

\noindent{\bf Proposition 5.3.3.} {\it
Let $G$ be a finite group and $p$ be the smallest prime dividing $|G|$.  If $|G'| \leq p^2 $, then every element of $G'$ is a commutator.
}

We conclude the thesis with a discussion, in the last chapter, on  some of the possible research problems related to  the results obtained in the earlier chapters.

\vspace{.5cm}

\end{document}